

\baselineskip=14pt
\parskip=10pt
\def\halmos{\hbox{\vrule height0.15cm width0.01cm\vbox{\hrule height
  0.01cm width0.2cm \vskip0.15cm \hrule height 0.01cm width0.2cm}\vrule
  height0.15cm width 0.01cm}}

\magnification=\magstephalf

\def\1{{\overline{1}}}
\def\2{{\overline{2}}}
\parindent=0pt
\overfullrule=0in

\def\frac#1#2{{#1 \over #2}}
\centerline
{\bf 
A combinatorial proof of Cramer's Rule
}
\bigskip
\centerline
{\it Doron ZEILBERGER}
\bigskip
\qquad {\it Dedicated to my Rutgers colleague Antoni A. Kosinski (b. May 25, 1930) who proved that Cramer's Rule is indeed due to Cramer. }

{\bf Prerequisites}: We assume that readers are familiar with the notions of : {\it integer} (e.g. $3$, $11$), 
{\it set} (e.g. $\{1,2,3\}$), {\it addition} (of numbers or symbols) (e.g. $a+b$) , {\it multiplication} (e.g. $bc$), and ``$> \, , \,  = \, , \, < \, , \, \leq \, , \, \geq$''.
Two symbols $a$ and $b$ {\it commute} if $ab=ba$. All the symbols in this paper commute.
No other knowledge is needed!

{\bf Notation}: If $w$ is a {\it weight} defined on a set $S$ (i.e. a way of assigning a number, or algebraic expression to members of $S$), then
$w(S):=\sum_{s \in S} w(s)$. For example if $S=\{1,4\}$, $w(1)=a$, $w(4)=b$, then $w(\{1,4\})=a+b$.

{\bf Definitions}

$\bullet$ For any positive integer $n$, an $n$-permutation is a list of integers $\pi=\pi_1 \dots \pi_n$,  where $1 \leq \pi_i \leq n$ for all $1 \leq i \leq n$
and $\pi_i \neq \pi_j$  if $i \neq j$. 

$\bullet$ The set of all $n$-permutations is denoted by $S_n$. For example, $S_3=\{123,132,213,231,312,321\}$.

$\bullet$ For an $n$-permutation $\pi=\pi_1 \dots \pi_n$, a pair $(i,j)$, where $1 \leq i<j \leq n$, is an {\it inversion} if $\pi_i>\pi_j$.
For example, If $\pi=51423$ then $(1,2),(1,3),(1,4), (1,5), (3,4), (3,5)$ are inversions of $\pi$. Let $inv(\pi)$ be the 
{\it number of inversions}, for example $inv(51423)=6$.

$\bullet$ Let $n$ be a positive integer, and  let $a_{i,j} (1 \leq i,j \leq n)$ and  $b_i$ ($1 \leq i \leq n$)
be $n^2+n$ {\it commuting} symbols (or numbers). Define $n+1$ {\it weights} $w_j$ ($0 \leq j \leq n)$ on $S_n$ as follows:
$$
w_0(\pi):= (-1)^{inv(\pi)} \prod_{1 \leq k \leq n } a_{\pi_k, k} \quad ,
$$
$$
w_j(\pi):= (-1)^{inv(\pi)} \, b_{\pi_j} \, \prod_{ {{1 \leq k \leq n} \atop {k \neq j}} } a_{\pi_k, k} \quad  (1 \leq j \leq n) \quad .
$$

{\bf Theorem} (Cramer [C]):  Let
$$
X_j :=w_j (S_n) \quad (0 \leq j \leq n) \quad,
$$
then
$$
x_j := \frac{X_j}{X_0} \quad (1 \leq j \leq n) \quad ,
$$
satisfy the $n$ linear equations
$$
\sum_{j=1}^{n} a_{i,j} x_j \, = \, b_i  \quad  (1 \leq i \leq n) \quad.
\eqno(C_i)
$$

{\bf Combinatorial Proof}:  By multiplying $(C_i)$ by $X_0$, we have to prove

$$
\sum_{j=1}^{n} a_{i,j} X_j \, = \, b_i X_0 \quad (1 \leq i \leq n) \quad.
\eqno(C'_i)
$$

For any positive integer $n$, let $F_n$ be the following set (with $n\cdot n!$ members, where $n!:=1 \cdot 2 \cdots n$)
$$
F_n:=\{[j, \pi] \, ; \, 1 \leq j \leq n \quad, \quad \pi \in S_ n \} \quad .
$$

For any $i$ ($1 \leq i \leq n$), define a weight $W_i$ on $F_n$ as follows:
$$
W_i([j,\pi]):=a_{i,j} w_j(\pi)=a_{i,j} (-1)^{inv(\pi)}  b_{\pi_j} \prod_{ {{1 \leq k \leq n} \atop {k \neq j}} } a_{\pi_k, k} \quad  .
$$

The left side of $(C'_i)$ is $W_i(F_n)$. 

{\bf Definition}: $[j,\pi] \in F_n$ is an  $i$-{\it good guy} if $\pi_j=i$, otherwise it is an {\it $i$-bad guy}. 
Let $G_{n,i}$ and $B_{n,i}$ be the subsets of $F_n$ consisting of the $i$-good guys and $i$-bad guys respectively.

Obviously, since $F_n=G_{n,i} \cup B_{n,i}$, we have
$$
W_i(F_n) \, = \, W_i(G_{n,i}) + W_i(B_{n,i}) \quad .
$$

{\bf Fact 1}: 
$$
W_i(G_{n,i})= b_i \, X_0 \quad
$$
{\bf Proof of Fact 1}:
If $[j,\pi] \in F_n$  is a good guy then, since $\pi_j=i$, we have:
$$
W_i([j,\pi])=(-1)^{inv(\pi)}  \, a_{i,j}  \,  b_{\pi_j}  \, \prod_{ {{1 \leq k \leq n} \atop {k \neq j}} } a_{\pi_k, k}
$$
$$
=(-1)^{inv(\pi)}  \, a_{\pi_j,j} \,  b_{i} \, \prod_{ {{1 \leq k \leq n} \atop {k \neq j}} } a_{\pi_k, k}
=(-1)^{inv(\pi)} \, b_{i}  \, \prod_{1 \leq k \leq n} a_{\pi_k, k}
\,= \, b_i w_0(\pi) \quad .
$$

Hence $W_i(G_{n,i})=b_i w_0(S_n)=b_i X_0$. \halmos

{\bf Fact 2}: 
$$
W_i(B_{n,i})= 0 \quad .
$$
{\bf Proof of Fact 2}: Let $[j,\pi] \in F_n$ be an $i$-bad guy. Let $a:=\pi_j$  and $j':=\pi^{-1}(i)$. Of course  $a \neq i$ and $j' \neq j$.
Define a permutation $\sigma$ by transposing
$\pi_j=a$ and $\pi_{j'}=i$, in other words
$\sigma_j=i$, $\sigma_{j'}=a$ and $\sigma_k=\pi_k$ if $k \not \in  \{j,j'\}$.
Let
$$
T_i([j,\pi]):=[j', \sigma] \quad .
$$
We have
$$
W_i([j,\pi]) \, = \, (-1)^{inv(\pi)} \, a_{i,j} \, b_a \, \prod_{{{1 \leq k \leq n} \atop { k \neq j}} } a_{\pi_k,k}
= (-1)^{inv(\pi)} a_{i,j} \, b_a \,  a_{\pi_{j'},j'} \, \prod_{{ {1 \leq k \leq n} \atop  {k \neq j,j'}}} a_{\pi_k,k}
$$
$$
= (-1)^{inv(\pi)} \, a_{i,j} \, b_a \, a_{i,j'} \, \prod_{ { {1 \leq k \leq n} \atop  {k \neq j,j'}} } a_{\pi_k,k} \quad .
$$
Similarly
$$
W_i([j',\sigma])= (-1)^{inv(\sigma)} \, a_{i,j'} \, b_a \, \prod_{{ {1 \leq k \leq n} \atop  {k \neq j'}} } a_{\sigma_k,k}
=(-1)^{inv(\sigma)} \, a_{i,j'} \, b_a \, a_{ \sigma_j,j} \, \prod_{ {{1 \leq k \leq n} \atop { k \neq j, j'} }} a_{\sigma_k,k}
$$
$$
=(-1)^{inv(\sigma)} \, a_{i,j'} \, b_a \, a_{i,j} \, \prod_{{ {1 \leq k \leq n} \atop  {k \neq j, j'} }} a_{\sigma_k,k} \quad .
$$

Since $inv(\pi)-inv(\sigma)$ is odd (why?), and $\pi_k$ and $\sigma_k$ coincide if $k \not \in  \{j,j'\}$,
we have, by {\it commutativity}, that for any $i$-bad guy $b$, $W_i(b)+W_i(T_i(b))=0$.

Since $T_i(T_i(b))=b$ for all  $i$-bad guys $b$ (why?), all the bad guys can be arranged into mutually $W_i$-canceling pairs, proving Fact 2. \halmos

Combining Facts 1 and 2,  $(C'_i)$, and hence Cramer's Rule, follow. \halmos

{\bf Comments}:  {\bf 1.} There are several `short' proofs of Cramer's rule that can be found in Wikipedia  and its references,
but they all assume  knowledge of linear algebra. Our proof is fully {\it self-contained}, and does not
assume {\it anything} besides the prerequisites listed at the beginning. We believe that if you include all the necessary background, our proof is the shortest.

{\bf 2.} For a fascinating defense of Gabriel Cramer's priority for his rule, see Antoni Kosinski's article [K].

{\bf References}

[C] Gabriel Cramer,  {\it ``Introduction à l'Analyse des lignes Courbes alg\'ebriques"}  Geneva (1750). pp. 656-659. \hfill \break
{\tt https://sites.math.rutgers.edu/\~{}zeilberg/mamarim/mamarimhtml/cramerCramer1750.pdf} \quad .

[K] A. A.  Kosinski, {\it  Cramer's Rule is due to Cramer}, Mathematics Magazine {\bf 74} (2001), 310-312. \hfill \break
{\tt https://sites.math.rutgers.edu/\~{}zeilberg/mamarim/mamarimhtml/cramerKosinski2001.pdf} \quad .

\bigskip
\hrule
\bigskip

Doron Zeilberger, Department of Mathematics, Rutgers University (New Brunswick), Hill Center-Busch Campus, 110 Frelinghuysen Rd., Piscataway, NJ 08854-8019, USA. \hfill\break
Email:  {\tt DoronZeil at gmail  dot com}   \quad .
\bigskip
Aug. 18, 2024.
\bigskip
{\bf Exclusively published in the Personal Journal of Shalosh B. Ekhad and Doron Zeilberger and arxiv.org}
\end